\documentstyle[12pt]{amsart}
\newtheorem{theo}{THEOREM}[section]

\newtheorem{prop}[theo]{Proposition}
\newtheorem{definition}[theo]{Definition}
\newtheorem{rem}[theo]{Remark}
\newtheorem{claim}[theo]{Claim}

\newtheorem{ex}[theo]{Example}
\newtheorem{facts}[theo]{FACTS}
\newcommand{\brref}[1]{(\ref{#1})}

\newcommand{\lra}{\longrightarrow}
\newcommand{\ra}{\rightarrow}
\newcommand{\emilia}{{Emilia MEZZETTI}}
\newcommand{\lucia}{{Maria Lucia FANIA}}
\newcommand{\e}{and}
\newenvironment{rem*}{\begin{rem}\em}{\end{rem}}
\newenvironment{ex*}{\begin{ex}\em}{\end{ex}}
\newenvironment{claim*}{\begin{claim}\em}{\end{claim}}
\newenvironment{definition*}{\begin{definition}\em}{\end{definition}}
\newenvironment{facts*}{\begin{facts}\em}{\end{facts}}

\title[]{On the Hilbert Scheme of Palatini threefolds}
\author[M.L. FANIA {\e} E. Mezzetti]{{\lucia$^*$} {\e} {\emilia$^*$}}
\address{ Maria Lucia Fania, Dipartimento di Matematica
\\Universit\`{a} degli Studi di L'Aquila\\
Via Vetoio Loc. Coppito\\67100 L'Aquila\\Italy}
\email{fania@@univaq.it}

\address{Emilia Mezzetti, Dipartimento di Scienze Matematiche
\\Universit\`{a} degli Studi di Trieste\\
Via Valerio 12/1\\
34127 Trieste\\
Italy\\}
\email{mezzette@@univ.trieste.it}
\thanks{ Partially supported by MURST in the framework
of the project:``Geometria sulle variet\`a algebriche"}

\begin{document}
\maketitle

%%%%%%%%%%%%% ABSTRACT %%%%%%%%%%%%%%%%%%%%%%%%%%%%%%%%%%%%%
\begin{abstract}

In this paper we study the Hilbert scheme of  Palatini threefolds $X$ in
${\Bbb P}^5$. We prove that such a scheme has an irreducible component
 containing $X$ which is birational to
the Grassmannian
${\Bbb G}(3,{\overset{\vee} {\Bbb P}}{^{14}})$ and we
 determine  the exceptional locus of the birational map.
\end{abstract}

%%%%%%%%%%%%%%% INTRODUCTION %%%%%%%%%%%%%%%%%%%%%%%%%%%%%

\section{introduction}

Let $X$ be a smooth non degenerate $3$-fold in ${\Bbb P}^5$
which is a scroll over a smooth surface $S$.
There are four examples, all classical, of such scrolls in ${\Bbb P}^5$:
the Segre scroll, the Bordiga scroll, the Palatini scroll and the $K3$-scroll,
of degree 3, 6, 7, 9 respectively.
Ottaviani in \cite{o}  proved that in ${\Bbb P}^5$ these are
the only smooth $3$-dimensional scrolls
over a surface.

The first two scrolls are arithmetically Cohen-Macaulay threefolds, defined
by the maximal minors
of a $3\times 2$, or of a $4\times 3$ matrix of linear forms respectively and
their Hilbert schemes are described by Ellingsrud (\cite{e}).

In this paper we are interested in studying  the Hilbert scheme
of Palatini scrolls.

By definition, a Palatini scroll is the degeneracy locus of a
general morphism $\phi:{\cal O}_{ {\Bbb P}^5}^{\oplus 4}\lra
{\Omega}_{{\Bbb P}^5}(2)$. This definition is the straightforward
analogous in ${\Bbb P}^5$ of the classical construction performed
by Guido Castelnuovo in 1891 for the projected Veronese surface
$S$ in ${\Bbb P}^4$ (see \cite{c}). In the modern language, he
showed indeed that any such surface $S$ can be interpreted as the
degeneracy locus of a morphism defined by three independent
sections of ${\Omega}_{{\Bbb P} ^4}(2)$.

This kind of morphisms are closely connected to the notion of
linear complexes of lines (see \cite{c}, \cite{bm}). In fact, it
is possible  to interpret the degeneracy locus of a general
morphism $\phi:{\cal O}_{ {\Bbb P}^5}^{\oplus 4}\lra
{\Omega}_{{\Bbb P}^5}(2)$ as the set of centres of linear
complexes in ${\Bbb P}^5$ belonging to a general linear system
$\Delta$ of dimension $3$ (web for short) of linear complexes in
${\Bbb P}^5$. This  latter interpretation will be very important
in proving our results.

In the case of Veronese surfaces, Castelnuovo proved also that the
net of complexes, giving raise to the morphism ${\cal O}_{ {\Bbb
P}^4}^{\oplus 3}\lra {\Omega}_{{\Bbb P}^4}(2)$, can be
reconstructed from $S$, as the span of the locus of trisecant
lines to $S$ inside the embedding space of the Grassmannian ${\Bbb
G}(1,4)$. This is equivalent to the claim that the irreducible
component of the Hilbert scheme of the Veronese surface is
birational to ${\Bbb G}(2,H^0({\Omega}_{{\Bbb P}^4}(2))$. Here we
study the analogous problem for Palatini threefolds.

One of the first things we prove, see Proposition \brref{Hilbertscheme of
Palatini}, is that the Hilbert scheme of the Palatini scroll has a
distinguished reduced irreducible component ${\cal H}$  which is smooth
at the point representing $X$ and of dimension $44$.  Such dimension is equal
to the dimension of ${\Bbb G}(3,{\overset{\vee} {\Bbb P}}{^{14}})$,
the Grassmannian parametrizing maps
${\cal O}_{ {\Bbb P}^5}^{\oplus
4}\lra {\Omega}_{{\Bbb P}^5}(2)$.

It's known that there is a natural map $\rho:{\Bbb
G}(3,{\overset{\vee} {\Bbb P}}{^{14}}) -- \ra {\cal H}$ and we
study such map. The main result we obtain is  the following
theorem. It relies on a careful description of the variety of
four-secant lines of the Palatini scroll. We prove that this variety 
is the union of $C_4$,  the base locus of a general
web of linear complexes, and of one more component, the one given by the
lines contained in $X$. Hence the situation is a bit different from the one
encountered by Castelnuovo in the case of a general net of linear complexes
in ${\Bbb P}^4$. But nevertheless we can reconstruct $C_4$ from $X$
 in a geometrical way, as it was done by Castelnuovo in the case of 
Veronese surfaces.

\begin{theo}
\label{birationa map}
Let $X \subset {\Bbb P}^5$ be a smooth  Palatini scroll of degree $7$. Let
${\cal H}$ be the irreducible component of the
Hilbert scheme containing $X$.
Then the rational map
$$\rho:{\Bbb G}(3,{\overset{\vee} {\Bbb P}}{^{14}}) -- \ra {\cal H}$$
 is birational.
\end{theo}

Since the  map $\rho$ is birational, our next task is to determine the
locus over which such map is not
regular.

It is well-known that if the degeneracy locus of a bundle map from a vector
bundle $F$ to a vector bundle $G$ has the
expected dimension then it lies in the same Hilbert
scheme as the degeneracy locus of a general map from $F$ to $G$. So $\rho$
is not regular over the points $\Delta
\in {\Bbb G}(3,{\overset{\vee} {\Bbb P}}{^{14}})$
such that
the set of centres of linear complexes in ${\Bbb P}^{5}$ belonging to $\Delta$
%corresponding degeneracy locus%
has dimension strictly bigger than $3$.

It turns out that these are the $3$-spaces of ${\overset{\vee} {\Bbb
P}}{^{14}}$
which either are completely contained in the dual Grassmannian
${\overset{\vee} {\Bbb G}}(1,5)$,
or intersect ${\Bbb G}(3,5)$, naturally identified with the singular locus
of ${\overset{\vee} {\Bbb G}}(1,5)$, along a curve. In other words, we can
consider
the intersection of $\Delta$ with ${\overset{\vee} {\Bbb G}}(1,5)$: in general
it is a cubic surface $S$ which can be identified with the base of the scroll.
The non-regularity of $\rho$ at $\Delta$ means that either $S$ is not defined
or that it is singular along a curve.

For the precise statements about the locus over which the map $\rho$ is not
regular, we  refer to \S 4,  Theorem \brref{caso A}
and Theorem \brref{caso B}.

In the last section we will see how the Hilbert scheme of the
Palatini scrolls fits into commutative diagrams involving the
variety of $6\times 6$ skew-symmetric matrices of linear forms,
the space of cubic surfaces, the moduli space of rank two bundles
$E$ on a cubic surface, with $c_1(E)={\cal O}_{S}(2), c_2(E)=5$.
The relations among the latter mathematical objects where
considered in \cite{bea} to which we refer for all the details.\\

{\it Aknowledgments.} The authors would like to thank the
referee for many interesting comments and suggestions.

\section{Notations and Preliminaries}

The following notation will be needed later on in the paper.
We refer to \cite{lb} for the details.

Let $X$ be a smooth non degenerate 3-fold in ${\Bbb P}^5$.
 By $Hilb^{4}{\Bbb P}^5$, respectively $Hilb^{4}X$, we denote the
Hilbert scheme of zero dimensional subschemes of ${\Bbb P}^5$,
respectively of $X$, of length $4$, i.e. $4$-tuples.

Let $Hilb^{4}_{c}{\Bbb P}^5$ be the open smooth subset of $Hilb^{4}{\Bbb P}^5$
of the 4-tuples lying on a smooth curve,
$Hilb^{4}_{c}X$:=$Hilb^{4}X \times_{Hilb^{4}P^{5}} Hilb^{4}_{c}{\Bbb P}^5$.
Note that dim $Hilb^{4}_{c}{\Bbb P}^5 = 20$.

Let $Al^{4}{\Bbb P}^5$ be the subvariety given by those elements in
$Hilb^{4}_{c}{\Bbb P}^5$,
which are on  some line in ${\Bbb P}^5$.

\begin{rem*}  (1) $Al^{4}{\Bbb P}^5$ is a smooth subvariety of
$Hilb^{4}_{c}{\Bbb P}^5$ of dimension 12.

(2) the map $a:Al^{4}{\Bbb P}^5 \lra {\Bbb G}(1,5)$ which sends an element of
$Al^{4}{\Bbb P}^5$ to the line on which it lies
is a fibration of fibre type $Hilb^{4}{\Bbb P}^1\cong {\Bbb P}^4$.
\end{rem*}

We recall the definition of the embedded $4$-secant variety of $X$.

Let $Al^{4}X$:=$Al^{4}{\Bbb P}^5 \times_{Hilb^{4}{\Bbb P}^5} Hilb^{4}_{c}X$.
We denote by $\Sigma_4(X)$ the closure of $a(Al^{4}X)$ in ${\Bbb G}(1,5)$. Let
$${\cal F}:=\{ (x, L)\in {\Bbb P}^5\times {\Bbb G}(1,5) | x\in L\}$$
be the flag manifold and let $p_1:{\cal F} \lra {\Bbb P}^5$,
$p_2:{\cal F} \lra {\Bbb G}(1,5)$ be
the two projections.

\begin{definition*}
\label{4secant variety}
$S_4(X):= p_{1}(p_{2}^{-1}(\Sigma_4(X)))\subset {\Bbb P}^5$ is the
embedded $4$-secant variety of $X$.
\end{definition*}

We recall a basic result on the Hilbert scheme which will be important
for us.

\noindent Let $Z$ be a smooth connected projective variety.
Let $X$ be a connected submanifold of $Z$ with $H^1(X,N)=0$ where
$N$ is the normal bundle of $X$. The following proposition holds,
see (\cite{gr}, \cite{so}).

\begin{prop}
\label{basic fact} Let $Z$ and $X$ be as above. There
exist irreducible projective varieties  ${\cal Y}$ and
${\cal H}$ with the following properties:
\begin{itemize}
\label{properties}
 \item[(i)]
${\cal Y} \subset {\cal H}\times Z$ and the map
$p:{\cal Y}\lra {\cal H}$ induced
by the product projection is a flat surjection,
 \item[(ii)]  there is a smooth point $x\in {\cal H}$
with $p$ of maximal rank in a neighborhood of $p^{-1}(x)$,
 \item[(iii)]  $q$ identifies  $p^{-1}(x)$ with $X$ where
$q:{\cal Y}\lra Z$ is the map induced
by the product projection, and
\item[(iv)] $H^0(N)$ is naturally identified with $T_{{\cal H},x}$ where
$T_{{\cal H},x}$ is the Zariski tangent space of ${\cal H}$ at $x$.
\end{itemize}
\end{prop}

%%%%%%%%%%%%%%%%%%%%%%%%%%%%%%%%%%%%%%%%%%%%%%%%%
%
%    Scroll di Palatini
%
%%%%%%%%%%%%%%%%%%%%%%%%%%%%%%%%%%%%%%%%%%%%%%%%%

\section{Webs of linear complexes in ${\Bbb P}^5$}

In this section we will study  the maps  ${\cal O}_{ {\Bbb P}^5}^{\oplus
4}\lra
{\Omega}_{{\Bbb P}^5}(2)$ and their degeneracy loci.

Note that a general morphism  $\phi:{\cal O}_{ {\Bbb P}^5}^{\oplus 4}\lra
{\Omega}_{{\Bbb P}^5}(2)$ is assigned by given $4$ general global sections of
${\Omega}_{{\Bbb P}^5}(2)$.
\vspace{3mm}

We  recall first two interpretations of global sections of
${\Omega}_{{\Bbb P}^5}(2)$, that explain the link with the
classical construction of Castelnuovo (\cite{c}).
 \vspace{3mm}

$\bullet$ Geometric interpretation of $\phi$ (see \cite{o}).
\vspace{3mm}

Let ${\Bbb P}^5={\Bbb P}(V)$. Consider the twisted dual Euler sequence on
${\Bbb P}^5$
\begin{eqnarray*}
0\ra {\Omega}_{{\Bbb P}^5}(2) \ra
\overset {6}{\oplus}{\cal O}_{{\Bbb P}^5}(1)
\ra {\cal O}_{{\Bbb P}^5}(2) \ra 0
\end{eqnarray*}

Being $V^*\cong  H^{0}({\Bbb P}^5,{\cal O}_{{\Bbb P}^5}(1))$, we have
the  natural identification
 Hom$({\cal O}_{{\Bbb P}^5},{\Omega}_{{\Bbb P}^5}(2))\cong
\wedge_{}^2V^*$,
where a morphism ${\cal O}_{{\Bbb P}^5}\ra {\Omega}_{{\Bbb P}^5}(2)$ is
given
(after choosing a basis in $V$ and its dual basis in $V^*$)
by a skew-symmetric $6\times 6$ matrix $A=(a_{ij})$,
$a_{ij}\in \Bbb C$
and corresponds to the morphism $V\ra V^*$ given by

$$(x_0,...,x_5) \ra (\sum a_{0i}x_i,...,\sum a_{5i}x_i)$$

The morphism $\phi:{\cal O}_{{\Bbb P}^5}^{4} \lra {\Omega}_{{\Bbb
P}^5}(2)$ is
given by four generic $6\times 6$ skew-symmetric
matrices $A, B, C, D$.

The equations of the degeneracy locus $X$ of $\phi$ are the $4\times 4$ minors
of the $4\times 6$ matrix
\begin{eqnarray}
\label{equation2}
F=
\left( \begin{array}{cccc}
\sum a_{0i}x_i&.&.&\sum a_{5i}x_i  \\
\sum b_{0i}x_i&.&.&\sum b_{5i}x_i \\
\sum c_{0i}x_i&.&.&\sum c_{5i}x_i\\
\sum d_{0i}x_i&.&.&\sum d_{5i}x_i
\end{array}
\right)
 \end{eqnarray}

If $P=(x_0,\ldots,x_5)\in X$ then there exists
$(x,y,z,t)\neq (0,0,0,0)$ such that

\begin{eqnarray}
\label{equation}
(xA+yB+zC+tD)P=0
\end{eqnarray}
or equivalently
\begin{eqnarray}
\label{equation1}
\sum_{i=0}^{5}
(xa_{ji}+ya_{ji}+ za_{ji}+ta_{ji})x_i = 0, \qquad j=0,...,5
\end{eqnarray}

Hence the matrix
$xA+yB+zC+tD$ has to be a
degenerate  skew-symmetric matrix and its pfaffian has to vanish.
Let $(x,y,z,t)$ be homogeneous
coordinates in ${\Bbb P}^3$, then the vanishing of the pfaffian of
$xA+yB+zC+tD$ defines a
hypersurface $S$ of degree $3$ in ${\Bbb P}^3$.

For $\phi$ general and a fixed $(x,y,z,t)\in S$, the matrix $xA+yB+zC+tD$
has rank four, so we find in $X$  a line of solutions of the
equation
\brref{equation1} and thus $X$ is a scroll over $S$.  Such $X$ is a Palatini
scroll of degree $7$, see (\cite{o}, \cite{bm}).

Let $f:X\to S$ denote the scroll map.
If we fix a point $P\in X$, then $f(P)$
is the unique solution of the equation
\brref{equation}, interpreted as an equation in
$(x,y,z,t)$. The unicity of the solution
is equivalent to the fact that
the $4\times 6$ matrix $F$  in \brref{equation2} has rank $3$.

The morphism $f:X\to S$ is associated to the line bundle $K_X+2H$,
(\cite{o}, 3.3), where $H$ is the hyperplane divisor. Moreover $X\cong {\Bbb P}_
{S}(E)$ where
$E:=f_{*}H$ is a rank two vector bundle on $S$ with
$c_1(E)={\cal O}_{S}(2), c_2(E)=5$.
\vspace{3mm}

$\bullet$ Global sections of
${\Omega}_{{\Bbb P}^5}(2)$ as linear  complexes.
\vspace{2.5mm}

Consider ${\Bbb G}(1,5)$, the Grassmannian of lines
in ${\Bbb P}^5$, embedded in ${\Bbb P}(\wedge_{}^2V) \cong {\Bbb P}^{14}$
via the Pl\"ucker map. The dual space
${\Bbb P}(\wedge_{}^2V^*)\cong {\overset{\vee} {\Bbb P}}{^{14}}$
parametrizes hyperplane sections of ${\Bbb G}(1,5)$ or, in the old
terminology,
linear complexes in ${\Bbb P}^5$. A linear complex $\Gamma$ in
${\Bbb P}^5$ is
 represented by a linear equation: $\sum_{0\le i<j\le 5} a_{ij}p_{ij}$ in
 the Pl\"ucker coordinates $p_{ij}$. We associate to $\Gamma$ the
skew-symmetric
matrix $A=(a_{ij})$ of order $6$.
A point $P\in {\Bbb P}^5$ is called a centre of $\Gamma$ if all lines through
$P$
belong to $\Gamma$. The space ${\Bbb P}(Ker (A))$ is the set of centres of
$\Gamma$
 and it is called the singular space of $\Gamma$. Since we are in ${\Bbb P}^5$
 a general linear complex $\Gamma$ does not have any centre. In fact let $A$ be
the skew-symmetric matrix associated to $\Gamma$. Being $\Gamma$ general
 it follows that $rk A=6$. $\Gamma$ is said to be special if its singular
space is at least a line. The special complexes can be of first type or of
second type depending on whether they have a line or a ${\Bbb P}^3$ as a
singular space.

Note that a  special $\Gamma$ corresponds to a tangent hyperplane section of
${\Bbb G}(1,5)$, or equivalently to a point of ${\overset{\vee} {\Bbb
P}}{^{14}}$
lying in
${\overset{\vee} {\Bbb G}}{(1,5)}$, the dual variety of ${\Bbb G}(1,5)$.

{}From this it follows that it is possible to interpret the degeneracy locus
$X$ of a general morphism
$\phi:{\cal O}_{ {\Bbb P}^5}^{\oplus 4}\lra
{\Omega}_{{\Bbb P}^5}(2)$ as the set of centres of linear complexes
 belonging to a general linear system $\Delta$ of dimension
$3$ (web for short) of linear complexes in
${\Bbb P}^5$. Since  $\Delta$ is general, it does not
contain any special complex of the second type. Moreover
such  $\Delta$ is spanned by four linearly
independent complexes $\Gamma_1,..., \Gamma_4$ and it corresponds to a
linear subspace ${\Bbb P}^3$ in
${\overset{\vee} {\Bbb P}}{^{14}}$:
 it intersects
${\overset{\vee} {\Bbb G}}{(1,5)}$ (which is a hypersurface of degree three)
 along a cubic
surface which is disjoint from ${\Bbb G}(3,5)$. Its
points represent the special complexes of $\Delta$ and so the surface
can be identified with $S$ (see \cite{bm}).

Hence given a general web of linear complexes in ${\Bbb P}^5$ its
set of centres is the degeneracy locus
$X$ of a general morphism
$\phi:{\cal O}_{ {\Bbb P}^5}^{\oplus 4}\lra
{\Omega}_{{\Bbb P}^5}(2)$ which is a Palatini scroll of degree $7$.
\vspace{1.5mm}

We see next that the Hilbert scheme of the Palatini scroll has a
distinguished reduced irreducible component of dimension equal to $44$.
\vspace{1.5mm}

\begin{prop}
\label{Hilbertscheme of Palatini}
Let $X \subset {\Bbb P}^5$ be the Palatini scroll of degree $7$.
Then the  Hilbert scheme of $X$ has an irreducible component, ${\cal
H}$,  which is smooth at the point
representing $X$ and of dimension $44$.
\end{prop}
\begin{pf} We have seen that $X\cong {\Bbb P}_{S}(E)$ where $E$
is a rank two vector bundle on a smooth cubic surface $S$ with
$c_1(E)={\cal O}_{S}(2), c_2(E)=5$. Note that $E$ is stable with
respect to ${\cal O}_S(1)$ and therefore simple: in fact $f_*H=E$,
where $f$ is the adjunction mapping of $X$, hence
$H^0(S,E(-1))=H^0(X,-(K+H))=0$.

Let $N$ denote the normal bundle of $X$ in ${\Bbb P}^5$.
 We will show that $H^1(X,N)=0$. Let
\begin{eqnarray}
\label{eulersequ}
0\lra {\cal O}_{X} \lra {\cal O}_{X}(1)^{\oplus 6}
\lra T_{{ {\Bbb P}^5}|{X}} \lra 0
\end{eqnarray}
be the Euler sequence on ${\Bbb P}^5$
restricted to $X$. From the cohomology sequence associated to
\brref{eulersequ} we have $h^i(X,T_{{ {\Bbb P}^5}|{X}})=0$
for $i\ge 1$. Using the following exact sequence
\begin{eqnarray}
\label{tangentsequ}
0\lra T_{X} \lra T_{{ {\Bbb P}^5}|{X}}
\lra N \lra 0
\end{eqnarray}
and the fact that $h^i(X,T_{{ {\Bbb P}^5}|{X}})=0$ for $i\ge 1$ we get that
\begin{eqnarray}
\label{cohomnormal}
h^{i}(X,N) = h^{i+1}(X,T_{X}) \qquad {\text {for} \quad  i\ge 1.}
\end{eqnarray}
Hence in order to compute $h^1(X,N)$ will be enough to compute
$h^{2}(X,T_{X})$.

Let $f:{\Bbb P}(E)\lra S$ be the scroll map, where $E$ is a rank two
vector bundle on $S$ as above.  We have the following sequences:
\begin{eqnarray}
\label{sequencescroll}
0\lra {\Omega}_{X|S} \lra (f^{*}E)(-1)
\lra {\cal O}_{X} \lra 0
\end{eqnarray}

\begin{eqnarray}
\label{relative tgbdl}
0\lra T_{X|S} \lra T_{X}
\lra f^{*}T_{S} \lra 0
\end{eqnarray}
Dualizing \brref{sequencescroll} we get
\begin{eqnarray}
\label{sequencedualized}
0\lra {\cal O}_{X}\lra (f^{*}E^*)(1)
\lra T_{X|S} \lra 0
\end{eqnarray}
By the projection formula and the fact that
$R^{1}f_{*}{\cal O}_{X}(1)=0$ it follows that
the 1-st direct image
$R^{1}f_{*}((f^{*}E^*)(1))=
E^{*}\otimes R^{1}f_{*}{\cal O}_{X}(1)=0$. Hence by the Leray spectral
sequence it follows that
$$H^i(X, (f^{*}E^*)(1))\cong H^i(S, f_{*}(f^{*}E^*)(1))\cong
H^i(S, E^*\otimes E)$$
for $i\ge 0$. In particular we have that
$h^3(X, (f^{*}E^*)(1))=0$ being  $h^3(S, E^*\otimes E)=0$.
On the other hand
 $h^i(X, {\cal O}_{X})=0$ for $i>0$.
Hence the cohomology sequence associated to \brref{sequencedualized}
gives that $h^3(X,T_{X|S})=0$.
We also know that
$H^i(X,f^{*}T_{S})\cong H^i(S,T_{S})$, $i\ge 0$ and that $h^3(S,T_{S})=0$.
Thus $h^3(X,f^{*}T_{S})=0$  and by \brref{relative tgbdl} it follows that
$h^3(X,T_{X})=0$. This latter fact along with \brref{cohomnormal}
gives $h^2(X,N)=0$.

We now compute $H^i(S,T_{S})$.
Since $S$ is a smooth cubic surface in ${\Bbb P}^3$, the Euler sequence
on ${\Bbb P}^3$ restricted to $S$ along with
\begin{eqnarray*}
\label{tangentsequcubic}
0\lra T_{S} \lra T_{{ {\Bbb P}^3}|{S}}
\lra {\cal O}_{S}(3) \lra 0
\end{eqnarray*}
and the fact that
$h^i(S,{\cal O}_{S}(3))=0$ for $i>0$ give that $h^2(S,T_{S})=0$ and thus
$h^2(X,f^{*}T_{S})=0$.
It remains to prove that $h^2(X,T_{X|S})=0$. From the cohomology
sequence associated to \brref{sequencedualized} it follows that
$h^2(X,T_{X|S})=h^2(S, E^*\otimes E)$. Note that
$$h^2(S, E^*\otimes E)=h^2(S, {\cal End}(E))=h^2(S, {\cal O}_{S})+
h^2(S, {\cal End}_0(E))$$
The last equality follows from the fact that
${\cal End}(E) \cong {\cal O}_{S}\oplus {\cal End}_0(E)$, where
${\cal End}_0(E)$ is the bundle of the traceless endomorphisms of $E$, see
 (\cite{br}, pg. 121). Since $S$ is a smooth
cubic surface in ${\Bbb P}^3$ it follows that $h^2(S, {\cal O}_S)=0$.
Moreover $h^2(S, {\cal End}_0(E))=0$, see (\cite{bea}, proof of Lemma 7.7).
Hence the cohomology
sequence associated to \brref{relative tgbdl} gives that
$h^2(X,T_{X})=0$ and thus, by \brref{cohomnormal}, $h^1(X,N)=0$.
In order to compute the dimension of ${\cal H}$, by (\brref{basic fact},
 (iv)) we need to compute $h^0(X,N)$.
But $h^0(N)=\chi(N)$ and by the Hirzebruch-Riemann-Roch theorem we know that
\begin{eqnarray}
\label{chiN}
\qquad \qquad \quad \chi(N) = \frac{1}{6}(n_1^3-3n_1n_2+3n_3)+
\frac{1}{4}c_1(n_1^2-2n_2)+\frac{1}{12}(c_1^2+c_2)n_1+
r\chi({\cal O}_{X})
\end{eqnarray}
where $n_i=c_i(N), c_i=c_i(X)$ and $r=rk(N)=2$.

Note that $n_3=0$ since $rk(N)=2$. We compute the remaining Chern classes
of $N$  from the exact sequence
\brref{tangentsequ} and we get:
$$n_1=K+6H; \quad  n_2=15H^2+6HK+K^2-c_2$$
The numerical invariants of the Palatini scroll and its
Hilbert polynomial have been computed in \cite{o} and they are:
$$KH^2=-8; K^2H=7; K^3=-2; -Kc_2=24$$

$$\chi({\cal O}_{X}(t))=
\frac{7}{6}t^3+2t^2+\frac{11}{6}t+1$$
{}From these we get that $c_2H=15$. Plugging these in
\brref{chiN} we get that $\chi(N)=44$ and thus $h^0(N)=\chi(N)=44$.
\end{pf}

\

Let $\Delta$  be a web of linear complexes in ${\Bbb P}^5$.
$\Delta$ is spanned by 4 independent linear complexes in ${\Bbb P}^5$
and hence it corresponds to a linear ${\Bbb P}^3$ in
${\Bbb P}({\wedge_{}^2 {V}^{*}})\cong  {\overset{\vee} {\Bbb P}}{^{14}}$.
Thus the webs of linear complexes in ${\Bbb P}^5$ are parametrized by
${\Bbb G}(3,{\overset{\vee} {\Bbb P}}{^{14}})$, the Grassmannian of
${\Bbb P}^3$'s in
${\overset{\vee} {\Bbb P}}{^{14}}$. Therefore we can define a natural
 rational map
$$\rho:{\Bbb G}(3,{\overset{\vee} {\Bbb P}}{^{14}}) --\ra {\cal H}ilb(X).$$
The map $\rho$ sends a general web $\Delta$ to its singular set $X$.
By Propositions \brref{basic fact} and
 \brref{Hilbertscheme of Palatini}, the Hilbert scheme of $X$,
${\cal H}ilb(X)$, has an irreducible component ${\cal H}$  which
is smooth at the point representing $X$ and of dimension $44$. It
is well known that the image of $\rho$
 is dense in ${\cal H}$. We let
$\rho$ denote also the map $\rho:{\Bbb G}(3,{\overset{\vee} {\Bbb
P}}{^{14}}) --\ra {\cal H}$. Since domain and codomain of $\rho$
have the same dimension $44$ and $\rho$ is dominant,  a general
fibre of $\rho$ is finite. So it is natural to ask for the degree
of such fibre.

An answer to this question is one of the main results of this paper.

 \begin{theo}
\label{birationa map}
Let $X \subset {\Bbb P}^5$ be a smooth  Palatini scroll of degree $7$. Let
${\cal H}$ be the irreducible component of the
Hilbert scheme containing $X$.
Then the rational map
$$\rho:{\Bbb G}(3,{\overset{\vee} {\Bbb P}}{^{14}}) -- \ra {\cal H}$$
 is birational.
\end{theo}

\vspace{2mm}

The  proof of Theorem
\brref{birationa map} will
follow at once after we have proved the following Claims and Proposition.

Let us point out that to prove this theorem we have tried to adapt
Castelnuovo's
result, (\cite{c}, \S 7, 8) to our case. Unfortunately his whole argument
does not
extend and the difficulty lies in the fact that the locus of the
$4$-secants of $X$ in ${\Bbb G}(1,5)$ contains  besides $C_4$, also
other lines.  By $C_4$ we denote the base locus of a general web $\Delta$ of
linear complexes in  ${\Bbb P}^5$ of which $X$ is the singular set. In fact
we  prove, see Proposition \brref{4secant variety},  that there is only one more
component, the one given by the lines of the ruling of $X$.
\vspace{1.5mm}

Let $\Delta$ be a general web of linear complexes. Let
$\Gamma_1, \Gamma_2, \Gamma_3, \Gamma_4$ be four
linear complexes of ${\Bbb P}^5$ generating $\Delta$.
 Let $C_4\subset {\Bbb G}(1,5)$ denote the base of the web
$\Delta$, that is, $C_4$ is the  family of lines in ${\Bbb P}^5$
which are common to $\Gamma_1, \Gamma_2, \Gamma_3, \Gamma_4$: it
is irreducible and $4$-dimensional, being the intersection of the
Grassmannian with a general linear space of codimension $4$. Let
$X$ be the set of centres of complexes belonging to $\Delta$. The
following claim is the analogous of the Castelnuovo's result about
trisecants of the Veronese surface in ${\Bbb P^4}$.

\begin{claim*}
\label{4secant}
The lines of $C_4$ are 4-secants of $X$.
\end{claim*}
\begin{pf}(of Claim)  The equations of the degeneracy locus
$X$ are the fifteen
$4\times 4$ minors  $F_{ij}$ of the $4\times 6$ matrix $F$, with $F$
as in \brref{equation2}.
If $F_{12}$ denotes the $4\times 4$ minor obtained by deleting the
1st and 2nd column of $F$ then
$$\{F_{12}=0\}$$
is the variety of degree $4$ which is made out of the lines of
$C_4$ which intersect the 3-plane $x_0=x_1=0$.

Let ${\it l}$ be a line in  $C_4$ which does not intersect the 3-plane
$x_0=x_1=0$.
Let  ${\it l} \cap \{F_{1 2}=0\}=\{P_1, P_2, P_3, P_4\}$:
 then each $P_i$ is also on
a line  ${\it r}$ of $C_4$ which intersects the $3$-plane $x_0=x_1=0$.
Hence there is a pencil of lines of $C_4$ through $P_i$: the pencil spanned by
${\it l}$ and ${\it r}$, and so $P_i\in X$.
This is true for every $P_i$, $i=1,...,4$, hence ${\it l}$ is a  4-secant
of $X$.
\end{pf}

\begin{claim*}
\label{orderone}
$C_4$ is a congruence of lines of order one, i.e. for a generic
point $x\in {\Bbb P}^5$,
 there is only one
line of $C_4$ passing through $x$.
\end{claim*}
\begin{pf}(of Claim)
In fact, as a cycle, $C_4$ coincides
with
$\sigma_1^4$, where $\sigma_1$ is the Schubert cycle of lines of ${\Bbb P}^5$
intersecting a fixed ${\Bbb P}^3$. Using Pieri formula, it is easy to show that
$\sigma_1^4=a_0\sigma_{4}+a_1\sigma_{3 1}+a_2\sigma_{2 2}=
\sigma_{4}+3\sigma_{3 1}+2\sigma_{2 2}$ (see \cite{dep}). Since the
coefficient
of $\sigma_{4}$ is equal to the number of lines of $\sigma_1^4$ passing
through
a general point, we have the claim.
\end{pf}

\begin{claim*}
\label{dim4secants}
$C_4$ is an irreducible component of  ${\Sigma}_4(X)$, the locus of all
$4$-secant lines of $X$ in the
Grassmannian ${\Bbb G}(1,5)$.
\end{claim*}

\begin{pf}(of Claim)
 Let  ${\Sigma}_{k}(X)$ denote the locus of all $k$-secant
lines of $X$ in the
Grassmannian ${\Bbb G}(1,5)$.
Note that every irreducible component of ${\Sigma}_4(X)$ has dimension $\le
5$.
If otherwise
then,  since
${\Sigma}_4(X)\subseteq {\Sigma}_2(X)$ and since ${\Sigma}_2(X)$
is irreducible of dimension $6$, it would follow that
 every secant line is also a $4$-secant,
which is impossible.

So either $C_4$ is a whole irreducible component of  ${\Sigma}_4(X)$,
or it is contained in an irreducible component $\Sigma '$ of dimension $5$.
But, in the second case, it follows from \cite{m}, Theorem 2.3, that
the lines of $\Sigma '$ cannot fill up ${\Bbb P}^5$, against the fact
that $C_4$ has order one.
\end{pf}

\begin{prop} 
\label{4secant variety}
$\Sigma_4(X)=C_4\cup\Sigma_{\infty}(X)$, where
$\Sigma_{\infty}(X)$ is
the variety of lines contained in $X$.
\end{prop}
\begin{pf} (of Proposition)
The number $q_4(X)$ of the $4$-secant lines  of $X$  passing
through a general point of ${\Bbb P}^5$ is finite by \cite{m}. Hence
$q_4(X)$ can be computed using the formula given by Kwak in \cite{kw} and it
turns out that
$q_4(X)=1$. So $C_4$ is the unique irreducible component of $\Sigma_4(X)$
whose lines fill up
${\Bbb P}^5$.  We have to exclude
the existence of another irreducible component $\Sigma '$, of dimension $4$
of $5$,
such that the union of the lines of $\Sigma'$ is strictly contained in
${\Bbb P}^5$.

If such a component  $\Sigma '$ exists, it follows from \cite{m} that $X$ 
contains
either a one-dimensional family of surfaces of ${\Bbb P}^3$ of degree at
least $4$, or
a $2$-dimensional family of plane curves of degree at least $4$.
In the first case, every hyperplane $H$ containing such a surface $S$
of ${\Bbb P}^3$ cuts
$X$ along a reducible surface: $H\cap X=S\cup S'$, where $S'$ is a surface
of degree
$7-k\leq 3$. So on $X$ there should be also a $2$-dimensional family
of surfaces of degree $\leq 3$. But then a general hyperplane section of
$X$ should contain
a $2$-dimensional family of curves of degree $\leq 3$.
This surface is a rational non-special surface of ${\Bbb P}^4$, which
has been extensively studied (see for instance \cite{A}): it can be easily
excluded that it contains such a family of curves.

In the second case, every $3$-space containing a plane curve on $X$,
of degree $k\geq 4$, should cut  $X$ residually in a curve of degree $\leq 3$.
So we would get on $X$ also a $3$-dimensional family of curves of degree
$\leq 3$.
Now, $X$ contains in fact both a three-dimensional family of conics and of
cubics,
but a computation in the Picard group of a general
hyperplane section of $X$ shows that their residual curves cannot be plane.

So both possibilities are excluded and the proof of the Proposition is
accomplished.
\end{pf}

\begin{pf}(of Theorem)
Let  $X=\rho(\Delta)$ be a general threefold in the image of $\rho$.
Then $C_4$ is the unique irreducible component of $\Sigma_4(X)$
whose lines fill up
${\Bbb P}^5$. This implies that $X$ comes via $\rho$ from a unique web
$\Delta$.
Therefore, by the theorem on the dimension of the fibres, there is an open
subset in ${\cal H}$ such that the fibres of $\rho$ over this open subset
are finite and of degree one.
This proves the theorem.
\end{pf}
\

%%%%%%%%%%%%%%%%%%%%%%%%%%%%%%%%%%%%%%%%%%%%%%%%%%%%%%%%%%%%%%%%%%%%%%%
%%  CAPIRE dove la mappa rho non e' regolare %%
%%%%%%%%%%%%%%%%%%%%%%%%%%%%%%%%%%%%%%%%%%%%%%%%%%%%%%%%%%%%%%%%%%%%%%%%

\section{Regularity of the map $\rho$}

 We will see that there are webs $\Delta$ over which the map $\rho$ is not
regular.
Our next task is to determine such webs.
%%Il fatto seguente segue dalla risoluzione data dal complesso di
%%%%%%Eagon-Northcott

It is well-known that if the degeneracy locus of a bundle map from a vector
bundle $F$ to a vector bundle $G$ has the
expected dimension then it lies in the same Hilbert
scheme as the degeneracy locus of a general map from $F$ to $G$. So $\rho$
is not regular over the webs $\Delta$ such that
the corresponding degeneracy locus has dimension strictly bigger than $3$.

 We recall the following facts about dual Grassmannians that  will be used in 
the sequel (see for
instance
\cite{ha}):
\begin{facts*}
\
\begin{itemize}
\label{fattisuGr}
 \item[$\bullet$ \ ]   ${\Bbb G}{(3,5)}$ can
 be naturally embedded in ${\overset{\vee} {\Bbb G}}$, where
 ${\overset{\vee} {\Bbb G}}$ stands for
the dual Grassmannian  ${\overset{\vee} {\Bbb G}}{(1,5)}$.
 We can interpret ${\Bbb G}{(3,5)}$ as the
set of singular complexes of the second type,  because a complex of second
type is determined uniquely by its singular space ${\Bbb P}^3$: it is formed
by the lines intersecting that ${\Bbb P}^3$.

\item[$\bullet$ \ ]   ${\overset{\vee} {\Bbb G}} = Sec({\Bbb G}{(3,5)})$,
the  variety of secant lines, and
${\Bbb G}{(3,5)}=Sing({\overset{\vee} {\Bbb G}})$, hence
${\overset{\vee} {\Bbb G}}=Sec(Sing({\overset{\vee} {\Bbb G}}))$.
\item[$\bullet$ \ ]   The linear spaces contained in ${\Bbb G}{(3,5)}$ have
dimension
$\le 4$. In particular, a linear ${\Bbb P}^3$ in ${\Bbb G}{(3,5)}$
represents the set of
the 3-spaces of a fixed ${\Bbb P}^4$ passing through a fixed point.
\end{itemize}
\end{facts*}

We will use Pl\"ucker coordinates $p_{ij}$ on ${\Bbb P}^{14}$ and the
dual coordinates $m_{ij}$ on its dual space.

Let $\Delta$ be a web of linear complexes
in ${\Bbb P}^5$.  Let  $\Gamma_1,..., \Gamma_4$ be four linearly
independent complexes which span $\Delta$. Hence $\Delta$ corresponds to
a ${\Bbb P}^3\subset {\overset{\vee} {\Bbb P}}{^{14}}$.
The special complexes of $\Delta$
are parametrized by  $\Delta\cap
{\overset{\vee} {\Bbb G}}{(1,5)}$.
In fact the space ${\overset{\vee} {\Bbb P}}{^{14}}$ parametrizes all linear
complexes: special complexes correspond to tangent hyperplane sections of
${{\Bbb G}}{(1,5)}$, that is to
points of ${\overset{\vee} {\Bbb G}}{(1,5)}$ which is a cubic hypersurface in
${\overset{\vee} {\Bbb P}}{^{14}}$.
Moreover special complexes of second type can be interpreted as
points  in ${{\Bbb G}}{(3,5)}$ (which is also embedded in
${\overset{\vee} {\Bbb P}}{^{14}}$), because a special complex of second type
is uniquely determined  by its singular space ${\Bbb P}^3$:
it is formed by the lines intersecting that ${\Bbb P}^3$.

Hence the following situations can occur:

\begin{itemize}
\label{possibili casi}
 \item[$\alpha$)] \ $\Delta\subset {\overset{\vee} {\Bbb G}}{(1,5)}$, that
is,
all the complexes of $\Delta$ are special, or
\item[$\beta$)]  \ $\Delta \not\subset{\overset{\vee} {\Bbb G}}{(1,5)}$ and
 $\Delta
\cap{\overset{\vee} {\Bbb G}}{(1,5)}$ is a cubic surface $S$, possibly
singular.
\end{itemize}

Let us consider case $\alpha$):
$\Delta\subset {\overset{\vee} {\Bbb G}}{(1,5)}$, that is, the case in
which all
the complexes of $\Delta$ are special.
Let $A,B,C,D$ be four $6\times 6$ skew-symmetric matrices associated
to the complexes $\Gamma_1,..., \Gamma_4$, respectively. Hence these matrices
 span $\Delta$.
Note that in this case for all
$(x,y,z,t)\in {\Bbb P}^3$ we have that $pf(xA+yB+zC+tD)\equiv 0$.
Hence  the equation $pf(xA+yB+zC+tD)= 0$ does not define any surface in
${\Bbb P}^3$.

\

Before stating our result concerning case $\alpha$)  we recall few facts,
already introduced in
\cite{bm}, which will be needed. Let
$$\psi: {\overset{\vee} {\Bbb G}}{(1,5)}-- \ra {\Bbb G}(1,5)$$
be the rational surjective map which
sends a special complex of first type to its singular line, see \cite{bm}.
So $\psi$ is regular on
${\overset{\vee} {\Bbb G}}{(1,5)} \setminus {\Bbb G}(3,5)$.
The  fibre $\psi^{-1}({\it l})$ is formed by the special complexes having
{\it l} as singular line. The closures of
these fibres are 5-dimensional linear spaces which are denoted by
 ${\Bbb P}_{\it l}^{5}$. In fact we may think of $\psi^{-1}({\it l})$ as
the linear system of hyperplanes in ${\Bbb P}^{14}$ containing the
tangent space to ${\Bbb G}(1,5)$ at the point {\it l}\ :
$T_{{\it l},{\Bbb G}(1,5)}\cong {\Bbb P}^8$,
see (\cite{bm}, section 3) for the details.

\begin{rem*} Let ${\it l}, {\it m}\in {\Bbb G}(1,5)$ be lines of ${\Bbb P}^5$.
Then the intersection of ${\Bbb P}_{\it l}^{5}$ with ${\Bbb G}(3,5)$ is a
smooth quadric of dimension $4$.
Also ${\Bbb P}_{\it l}^{5}\cap {\Bbb P}_{\it m}^{5}$ is contained in
${\Bbb G}(3,5)$ and is just one point if {\it l}, {\it m} do not intersect,
or a plane if they intersect, see (\cite{bm}, Remark 1).
\end{rem*}

\begin{theo}
\label{caso A}
Let $\Delta$ be a web of linear complexes in ${\Bbb P}^5$ as
in case $\alpha$).
Then the rational map
$$\rho:{\Bbb G}(3,{\overset{\vee} {\Bbb P}}{^{14}}) -- \ra {\cal H}$$
 is not regular at the point corresponding to $\Delta$.
\end{theo}
\begin{pf}
Since the web $\Delta$ is as in case $\alpha$) then the following situations
can occur:
\begin{itemize}
\label{2 sottocasi}
 \item[$\alpha1$) \ ]   $\Delta\subset {\Bbb G}{(3,5)}$, or
 \item[$\alpha2$) \ ]  $\Delta\not\subset {\Bbb G}{(3,5)}$ and
$\Delta\subset{\overset{\vee} {\Bbb G}}$.
\end{itemize}

Let $\Delta$ be as in case  $\alpha1$), i.e.  $\Delta\subset {\Bbb G}{(3,5)}$.
Then by \brref{fattisuGr} we see that $\Delta$  represents
the set of the 3-dimensional linear spaces of a fixed ${\Bbb P}^4$
passing through
a fixed point. Let us fix the flag: $E_0\subset H_5$, where
$E_0=(1,0,0,0,0,0)\in {\Bbb P}^5$
and $H_5$ is the hyperplane in ${\Bbb P}^5$ whose equation is $x_5=0$.

Hence $\Delta$
represents the 3-dimensional linear spaces of ${\Bbb P}^5$ whose equations are:
$x_5=0, xx_1+yx_2+zx_3+tx_4=0$, with $(x,y,z,t)\in {\Bbb P}^3$.
The lines intersecting such ${\Bbb P}^3$ satisfy the following equations
in the Pl\"ucker coordinates $p_{ij}$:
 $xp_{15}+yp_{25}+zp_{35}+tp_{45}=0$. Thus $\Delta$ has
equations:
$m_{01}=m_{02}=m_{03}=m_{04}=m_{05}=m_{12}=m_{13}=m_{14}=m_{23}=m_{24}=
m_{34}=0$ in ${\overset{\vee} {\Bbb P}}{^{14}}$.

The corresponding matrix is
$$\left( \begin{array}{cccccc}
0&0&0&0&0&0\\
0&0&0&0&0&x\\0&0&0&0&0&y\\0&0&0&0&0&z\\
0&0&0&0&0&t\\0&-x&-y&-z&-t&0
\end{array}
\right)$$

Hence a generic element of $\Delta$ can be written in the form
$xA+yB+zC+tD$, where
$(x,y,z,t)\in {\Bbb P}^3$ and $A,B,C,D$ are $6\times 6$ skew-symmetric
matrices. Note that
$$A=\left( \begin{array}{cccccc}
0&0&0&0&0&0\\
0&0&0&0&0&1\\0&0&0&0&0&0\\0&0&0&0&0&0\\
0&0&0&0&0&0\\0&-1&0&0&0&0
\end{array}
\right)$$
and similarly we can write down $B,C,D$.
So $X$ is the variety whose equations are the
$4\times 4$ minors  $F_{ij}$ of the $4\times 6$ matrix
\begin{eqnarray}
F=
\left( \begin{array}{cccc}
\sum a_{0i}x_i&.&.&\sum a_{5i}x_i  \\
\sum b_{0i}x_i&.&.&\sum b_{5i}x_i \\
\sum c_{0i}x_i&.&.&\sum c_{5i}x_i\\
\sum d_{0i}x_i&.&.&\sum d_{5i}x_i
\end{array}
\right) =
\left( \begin{array}{cccccc}
0&x_5&0&0&0&-x_1\\
0&0&x_5&0&0&-x_2 \\
0&0&0&x_5&0&-x_3\\
0&0&0&0&x_5&-x_4
\end{array}
\right)
 \end{eqnarray}

It is straightforward to see that $X=3H_5$, where $H_5: x_5=0$.
Hence the map $\rho$ is not regular in this case.

\
Let $\Delta$ be as in case  $\alpha2$), i.e.  $\Delta\not\subset {\Bbb
G}{(3,5)}$ and
$\Delta\subset{\overset{\vee} {\Bbb G}}$.
According to the position of
$\Delta$ with respect to the fibre of
$\psi:{\overset{\vee} {\Bbb G}}{(1,5)}--\ra {\Bbb G}(1,5)$, we have
the following  situations:

\begin{itemize}
\label{4 sottocasi}
 \item[$\alpha$2.1) \ ]   there exists a unique line {\it l} such that
$\Delta\subset {\Bbb P}_{\it l}^{5}$,
 \item[$\alpha$2.2) \ ] $\Delta\not\subset {\Bbb P}_{\it l}^{5}$, for every
{\it
l}, but
there exists a line {\it m} such that $\Delta\cap {\Bbb P}_{\it m}^{5}=\pi$,
 with $\pi$ a plane,
\item[$\alpha$2.3) \ ] the general fibre of $\psi_{|\Delta}$ is a line,
\item[$\alpha$2.4) \ ] the general fibre of $\psi_{|\Delta}$ is a point.
\end{itemize}
\vspace{2.5mm}

In the case $\alpha$2.1) since $\Delta\subset {\Bbb P}_{\it l}^{5}$, then
$\Delta\cap{\Bbb G}{(3,5)}=
\Delta\cap({\Bbb P}_{\it l}^{5}\cap{\Bbb G}{(3,5)})$. Note that
${\Bbb P}_{\it l}^{5}\cap{\Bbb G}{(3,5)}$ is a smooth 4-dimensional quadric,
$\Delta$ is a ${\Bbb P}^3$ and thus $\Delta\cap{\Bbb G}{(3,5)}$ is a
2-dimensional quadric of ${\Bbb P}^3$: its points correspond to complexes
having
a ${\Bbb P}^3$ containing ${\it l}$ as a singular set, these ${\Bbb P}^3$ are
contained in the degeneracy locus $X$ and thus dim$ X > 3$.
\vspace{2.5mm}

In the case $\alpha$2.2) since there exists a line {\it m} such that
$\Delta\cap {\Bbb P}_{\it m}^{5}=\pi$,
 with $\pi$ a plane, then
$\Delta\cap{\Bbb G}{(3,5)}\supset (\Delta\cap{\Bbb P}_{\it m}^{5})\cap
({\Bbb G}{(3,5)}\cap {\Bbb P}_{\it m}^{5})$. This latter intersection is either
$\pi$ or a conic, since ${\Bbb G}{(3,5)}\cap{\Bbb P}_{\it m}^{5}$ is a
smooth 4-dimensional quadric and $\Delta\cap {\Bbb P}_{\it m}^{5}=\pi$.

As in the previous case, the points in such intersection correspond
to complexes having a
${\Bbb P}^3$ as a singular set which is contained in the degeneracy locus $X$
 and thus dim$ X > 3$.
\vspace{2.5mm}

In the case $\alpha$2.3) the general fibre of
$\psi_{|\Delta}: \Delta --\ra \psi_{|\Delta}(\Delta)\subset{\Bbb G}(1,5)$ is
a line.
Thus $T=\psi_{|\Delta}(\Delta)$ is a surface in ${\Bbb G}(1,5)$.
The intersection
$(\Delta\cap{\Bbb P}_{\it l}^{5})\cap
({\Bbb G}{(3,5)}\cap {\Bbb P}_{\it l}^{5})\neq\emptyset$, for every
${\it l}\in T$.  Hence for every ${\it l}\in T$, we find a ${\Bbb P}^3$
 contained in $X$ and thus dim$ X > 3$.
\vspace{0.5mm}

In the case $\alpha$2.4) the general fibre of
$\psi_{|\Delta}: \Delta --\ra\psi_{|\Delta}(\Delta)\subset {\Bbb G}(1,5)$
is a point. This means that $T=\psi_{|\Delta}(\Delta)$ is a 3-fold in
${\Bbb G}(1,5)$ and that
a general complex in $\Delta$ has a line as singular set. Thus
there is a 3-dimensional family of such lines and the variety $X$
is their union.

It could a priori happen that
$\dim X=3$, but then infinitely many lines of the family pass through a
general point of
$X$. Hence the matrix  $F$ (see \brref{equation2})
should have rank $<3$ at every point of $X$. This means that the four
hyperplanes,
whose coordinates are the rows of $F$, belong to a pencil, whose support is a
$3$-space $\Lambda_P$: it is the union of the lines of $C_4$ (the base of the
web $\Delta$) passing through $P$.
Assume that the first two rows of $F$ are linearly independent, then the
Pl\"ucker
coordinates of
$\Lambda_P$ are the order two minors of the first two rows of $F$.
 If $F_{0123}$  is the  minor of the last two columns, then its vanishing
at $P$
represents the condition that
$\Lambda_P$ meets the line $x_0=x_1=x_2=x_3=0$ (as in Claim 3.4).
If $\Lambda_P$ is disjoint from this line, let
$Q_P$ be the quadric surface $\Lambda_P \cap\{F_{0123}=0\}$: we claim that
$Q_P\subset X$.
Indeed, if
$z\in Q_P$, then both the line of $C_4$ through $P$ and the $3$-plane
$\Lambda_z$ contain
$z$. Let $P'$ be the point $\Lambda_z\cap \{x_0=x_1=x_2=x_3=0\}$: the line
$\overline{zP'}$ is in $C_4$ but not in $\Lambda_P$, so $z\in X$, and the
claim is proved.

Let now $P$, $P'$ be two distinct points in $X$: if
$\Lambda_P=\Lambda_P'$, then
$\Lambda_P \subset X$: in fact if $r$ is a line of $\Lambda_P$
through $P$, then  $r \in C_4$, so if $R \in r$, then $R \in \Lambda_P'$
so the line $RP'\in C_4$ and $R \in X$. Therefore we can deduce that, if
$\dim X=3$,  to each point $P\in X$ we can associate a ${\Bbb P}^3$
$\Lambda_P$ and
different points give different $3$-spaces. So also the quadrics $Q_P$ are  two
by two distinct and $X$ contains a family of dimension $3$ of quadric
surfaces: a
contradiction. Hence we conclude that $\dim X>3$.
\end{pf}
\vspace{1mm}

Before considering  case $\beta$) we make the following remark
about equations of cubic surfaces in ${\Bbb P}^3$ which will be used
later on in the paper.
\vspace{1mm}

\begin{rem*}\label{beauv} Let $S$ be a cubic surface in ${\Bbb P}^3$.
If $S$ is smooth then it can be defined by an equation $pfM=0$, where $M$ is a
$6\times 6$ skew-symmetric matrix of linear forms (\cite{bea}, Prop. 7.6).

Let $S$ be a singular element of $|{\cal O}_{P^3}(3)|$.
If $S$ has a finite number of lines then  by
(\cite {bl}) it  follows that its equation can  be expressed as
the determinant of a $3\times 3$  matrix $N$ of linear forms, except if $S$
is the surface whose class of projective equivalence  is called
$T_1$. Note that if $S$
is defined by $detN=0$,
then $S$ is also defined by $pf M=0$, where
  $M= \left( \begin{array}{cc}
0&N\\-^{t}N&0
\end{array}
\right)$,
since $pf M= detN$. Up to automorphisms the equation of the surface
in the class $T_1$ is
the following:
$x_0x_1^2+x_1x_3^2+x_2^3=0$. It is easy to check that it is the pfaffian of
the following
matrix
\vspace{1mm}

$$ \left( \begin{array}{cccccc}
0&-x_0&0&0&x_2&x_1\\
x_0&0&-x_0&x_3&0&x_2\\0&x_0&0&x_2&x_3&0\\0&-x_3&-x_2&0&x_1&0\\
-x_2&0&-x_3&-x_1&0&0\\-x_1&-x_2&0&0&0&0
\end{array}
\right)$$
\vspace{1mm}

If $S$ has an infinite number of lines then by a well known result
(see for instance \cite{fc}) $S$ is either a reducible cubic surface,
or an irreducible cone, or an irreducible cubic surface with a double line
(i.e a general ruled cubic surface).
But in all of these cases $S$ can be defined by
an equation
$detN$, where $N$ is a $3\times 3$  matrix of linear forms, and
hence also by $pfM=0$, with $M$ as above.

Hence we can conclude that every cubic surface can be expressed
as a pfaffian.
\end{rem*}

Let us consider  now  case $\beta$). Let $S$ be the cubic surface $\Delta
\cap{\overset{\vee} {\Bbb G}}{(1,5)}$. One of the following may
happen:
\begin{itemize}
\label{3 sottocasi}
 \item[$\beta$1)]  \ $\Delta\cap {\Bbb G}{(3,5)}=\emptyset$, or
 \item[$\beta$2)] \ $\Delta\cap {\Bbb G}{(3,5)}\neq \emptyset$, but
$S\not\subset
{\Bbb G}{(3,5)}=Sing({\overset{\vee} {\Bbb G}})$, or
\item[$\beta$3)] \ $\Delta\cap {\Bbb G}{(3,5)}\neq \emptyset$ and
$S\subset {\Bbb G}{(3,5)}$, that is, $S\subset Sing({\overset{\vee} {\Bbb
G}})$.
\end{itemize}
\vspace{2mm}

Let $\Delta$ be as in $\beta$1). Then  either $S$ is smooth, or it is
singular and its singularities correspond to tangency points of
$\Delta$ to ${\overset{\vee}
{\Bbb G}}{(1,5)}$. The case of $S$ smooth was  considered in Theorem
\brref{birationa map}. We give next an explicit example of a surface $S=\Delta
\cap{\overset{\vee} {\Bbb G}}{(1,5)}$ with
$\Delta\cap {\Bbb G}{(3,5)}=\emptyset$ and $S$ singular.

\begin{ex*}
\label{es}
Let $S$ be the cubic surface whose equation is
$x^2y-x^2z-xy^2+xz^2+y^3-y^2t+yzt=0$. The
class of projective equivalence of this surface is called
$T_4$ in \cite{bl}.  Its equation
can be written as  $\det N=0$, where $N$ is the following matrix:
$ \left( \begin{array}{ccc}
t&x&y\\
y+z&-y&2x+t\\
y&0&x+y-z
\end{array}
\right)$
or as $pf M=0$, where $M$ is obtained from $N$ as in \brref{beauv}.
\vspace{1mm}

\noindent $S$ has only one
singularity at the point $(0,0,0,1)$, which corresponds to the only matrix
of rank
less than $6$ obtained from $M$, for a particular choice of $x,y,z,t$.
It is easy to see that its rank is $4$.
\end{ex*}

Note that if  $\Delta$ is as in $\beta$1), for all $(x,y,z,t)\in S$, the matrix
$xA+yB+zC+tD$ has rank
four, so it determines always a line of solutions of the equation
$(xA+yB+zC+tD)^t(x_0 \dots x_5)=0$ on the degeneracy locus
$X$. Thus $\dim X=3$.

\begin{rem*}\label{formedip}
It is possible  that, for certain surfaces $S$, the matrices $A,B,C,D,$
appearing in a
pfaffian $pf(xA+yB+zC+tD)$ giving the equation of $S$, are linearly
dependent, so they
don't generate a $3$-space $\Delta$ in ${\overset{\vee} {\Bbb
P}}{^{14}}$, but only a plane.

For example, let $S$ be the
union of $3$ planes, then $S$ is defined by
$det(M)=0$, where $M$ is a $3\times 3$ diagonal
matrix whose non zero entries are the linear forms defining the three planes:
$F=ax+by+cz+td, G=a'x+b'y+c'z+t'd, H=a''x+b''y+c''z+t''d$.

Note that $S$ is also defined by $pf (xA+yB+zC+tD)=0$, where
$$A=\left( \begin{array}{cccccc}
0&0&0&a&0&0\\
0&0&0&0&a'&0\\0&0&0&0&0&{a''}\\-a&0&0&0&0&0\\
0&-a'&0&0&0&0\\0&0&-a''&0&0&0
\end{array}
\right)$$

The remaining matrices $B,C,D$ are of the same type where the entries
$a,a',a''$ are replaced by
 $b,b',b''$,  $c,c',c''$,  $d,d',d''$,  respectively.
Thus the  matrices $A,B,C,D$ are linearly dependent.

Another example is that of cones: if $S$ is a cone of vertex $(0,0,0,1)$
over a smooth
plane elliptic curve, then its equation can be put in  Weierstrass normal form
$y^2z=x(x-z)(x-cz)$: this can be interpreted as the determinant of
a $3\times 3$
matrix whose entries are linear forms in the variables $x, y, z$ only.

Nevertheless, for both examples, it is possible to find also another
pfaffian expression of
the equation of $S$, with a matrix which is a linear combination of four
independent matrices.
In the first case, for example, the equation of $S$ is
the pfaffian of the matrix $M$:
$$\left( \begin{array}{cccccc}
0&0&0&F&0&0\\
0&0&0&0&G&x\\0&0&0&0&0&H\\-F&0&0&0&0&0\\
0&-G&0&0&0&0\\0&-x&-H&0&0&0
\end{array}
\right)$$
In the second case, the equation of $S$ is the pfaffian of
$$\left( \begin{array}{cccccc}
0&x&t&y&0&0\\
-x&0&0&0&y&y\\-t&0&0&x&-z&0\\-y&0&-x&0&-l&-l\\
0&-y&z&l&0&x\\0&-y&0&l&-x&0
\end{array}
\right)$$
where $l=(c+1)x-cz$.
\end{rem*}

We prove next that the case $\beta$3) does not occur.

\begin{prop}
\label{B3nonsuccede}
 The case $\beta$3) does not occur.
\end{prop}
\begin{pf} Let $S$ be as in $\beta$3). Let $P,Q\in S$ and let $L$ denote
the line
through $P$ and $Q$. Then
$L\subset \Delta$. Moreover $L\subset {\overset{\vee} {\Bbb G}}$ since
${\overset{\vee} {\Bbb G}}= Sec({\Bbb G}(3,5))$. Thus
$L\subset S$ and hence $Sec(S)=S$. This latter fact implies that (the
support of) $S$ is
linear, hence it is a triple plane. So $\Delta$ should be tangent to
${\Bbb G}(3,5)$
along a plane. But this is impossible, because every tangent space to
${\Bbb G}(3,5)$ is
tangent at one point only, and the tangent spaces intersect two by two
along a plane.
\end{pf}

\

Let us consider now the case $\beta$2), that is: $\Delta\cap {\Bbb
G}{(3,5)}\neq
\emptyset$ and
$S\cap {\Bbb G}{(3,5)}=S\cap Sing({\overset{\vee} {\Bbb G}})\neq \emptyset$,
where $S$ is the cubic surface $S=\Delta\cap {\overset{\vee} {\Bbb G}}$, which
in this case is singular.

For all points $(x,y,z,t)$ in  $\Delta\cap {\Bbb
G}{(3,5)}$, the equation $$(xA+yB+zC+tD) \ ^t(x_0
\dots x_5)=0$$ is satisfied by the points of a ${\Bbb P}^3$, which enters
in the degeneracy
locus $X$. If the  intersection $\Delta\cap {\Bbb
G}{(3,5)}$ is a finite set of points, say $d$ points, then $X$
has dimension $3$ and contains
$d$ $3$-spaces as irreducible components. This number $d$ is  at most
$4$
(see \cite{bl}).  If  the intersection is infinite, then it contains at
least a curve $C$,
and over every point of $C$ there is a $3$-space contained in $X$,
therefore $\dim X\geq 4$.

The latter situation can appear only if $S$ is irreducible with a
double line,
or reducible in
the union of a plane with a quadric (possibly reducible). We will
give examples of
both these situations.

\begin{ex*}\label{es2}

$(i)$ Let $S$ be a cubic surface  having the line $r$ of equation $x=y=0$
as double line, and
containing also the lines $y=z=0$ and $x=t=0$. Then its equation takes the
form $F(x,y,z,t)=0$, where
$$F(x,y,z,t)=ax^2y+bx^2z+cxy^2+dxyz+exyt+fy^2t$$
Note that $F=det\ N$, where
$ N=\left( \begin{array}{ccc}
ex+fy&bx+dy&ax+cy\\
0&-y&z\\
-x&0&t
\end{array}
\right)$,
or $F=pf M$ where $M$ is the corresponding skew-symmetric matrix
$xA+yB+zC+tD$ (as in \brref{beauv}).
It is easy to check that the rank of $M$ is less than $4$ precisely at
the points  of $r$. Therefore, if $\Delta$ is the $3$-space generated by
$A,B,C,D$, $r$ is the intersection $\Delta\cap{\Bbb G}{(3,5)}$.

\

$(ii)$ Let now $S$ be the union of the quadric  $Q:xz-yt=0$ with the plane
$\pi:L(x,y,z,t)=0$. Then the determinant $ \left| \begin{array}{ccc}
x&y&0\\
t&z&0\\
0&0&L
\end{array}
\right|$ clearly vanishes on $S$. Hence the rank of the corresponding
skew-symmetric matrix
$M$ is $2$ along the conic $Q\cap \pi$. If we replace $Q$ with a quadric
cone $Q'$, with
similar computations we get rank $2$ on the conic $Q'\cap\pi$ and moreover
in the vertex of
$Q'$.
\end{ex*}

We conclude this section collecting the results so far obtained about
case $\beta$) in the
following Theorem:

\begin{theo}
\label{caso B}
Let $\Delta$ be a web of linear line complexes in ${\Bbb P}^5$ as
in case $\beta$).
Then the rational map
$$\rho:{\Bbb G}(3,{\overset{\vee} {\Bbb P}}{^{14}}) -- \ra {\cal H}$$
 is not regular at the point corresponding to $\Delta$ if and only if the
intersection
of $\Delta$ with ${\Bbb G}(3,5)$ contains a line or a conic.
\end{theo}

%%%%%%%%%%%%%%%%%%%%%%%%%%%%%%%%%%%%%%%
%%%%%  mappa pf di Beauville  %%%%%%%%%
%%
%%%%%%%%%%%%%%%%%%%%%%%%%%%%%%%%%%%%%%%
\section{The pfaffian map}

Let ${\cal S}_3$ be the variety of the  $6\times 6$ skew-symmetric
matrices of linear forms on ${\Bbb P}^3$. Let $pf:{\cal
S}_{3}--\ra |{\cal O}_{{\Bbb P}^3}(3)|$ be the rational map which
sends $M\in{\cal S}_{3}$ to the cubic surface $S$ defined by $pf
M=0$. It is regular on the open set of matrices whose pfaffian is
not identically zero.
 Note that the equation $pf M=0$ is equivalent to $pf (^{t}PMP)=0$,  with
 $P\in GL(6)$.  In fact
 $pf (^{t}PMP)=(detP) (pf M)$ and thus $S$ can also be defined by
$pf (^{t}PMP)=0$.

Let $GL(6)$ act on ${\cal S}_3$ by congruence, that
is the action is given by $^{t}PMP$ where $M\in {\cal S}_{3}$ and $P\in
GL(6)$.
As noted in \cite{bea}
the group $GL(6)$ acts freely and properly on
 ${\cal S}_{3}$ and the map $pf$
factors through ${\cal S}_{3}/GL(6)$.  Hence we have the following
commutative diagram:
\

\

$\put(60,7){${\cal S}_3/GL(6)$}
%\put(114,10){\vector(1,0){40}}
\put(129,10){\line(1,0){10}}
\put(144,10){\vector(1,0){10}}
\put(114,10){\line(1,0){10}}
\put(120,-2){$\overline{pf}$}
\put(158,7) {$|{\cal O}_{{\Bbb P}^3}(3)|$}
%\put(170,39){\vector(0,-1){20}}
\put(170,39){\line(0,-1){8}}
\put(170,27){\vector(0,-1){8}}
\put(164,45){${\cal S}_3$}
\put(172,27){$pf$}
%\put(155,45){\vector(-3,-2){35}}
\put(155,45){\line(-3,-2){10}}
\put(140,35){\line(-3,-2){10}}
\put(125,25){\vector(-3,-2){10}}
\put(125,37){$\pi$}$

\

\noindent  It's easy to
see that dim ${\cal S}_{3}/GL(6) = 60 - 36 = 24$. By Remark \brref{beauv}
it follows
that
 the pfaffian map, and by consequence also $\overline{pf}$
 is surjective.

If $S$ is a generic cubic surface defined by an equation $pf M=0$,
 Beauville in \cite{bea} has also proved that expressing $S$ as a
 linear pfaffian is equivalent  to
produce a rank $2$ vector bundle $E_{M}$ on
 $S$ which is so defined
$$E_{M}:= coker({\cal O}_{{\Bbb P}^3}(-1)^6 \buildrel{M}\over{\lra}
 {\cal O}_{{\Bbb P}^3}^6).$$

Let  ${\cal S}_3^{s}$ be the scheme of matrices $M\in{\cal S}_3$
such that $pf M$ is smooth: it follows that the map $M$ is
injective and $E_M$ is a stable bundle. Then ${\cal S}_3^{s}$ is
an open non-empty set of ${\cal S}_3$. Going through the proof of
Lemma 7.7 in \cite{bea}, one gets that for the generic $S$ the
fibre $\overline{pf}^{-1}(S)$ can be identified with an open
subset of the moduli space of simple rank $2$ vector bundles $E$
on $S$ with $c_1(E)={\cal O}_{S}(2), c_2(E)=5$.
 These are precisely the
bundles on $S$
such that ${\Bbb P}(E)$ is a Palatini scroll over $S$ (\cite {o}).
%Note that the dimension of these fibres is $$\dim {\cal S}_3/GL(6)-\dim |{\cal
%O}_{{\Bbb P}^3}(3)|=24 - 19=5.$$\

 Moreover
the quotient ${\cal S}_3^{s}/GL(6)$ can be canonically identified
with the set of pairs $(S,E)$, where $S$ is a  smooth cubic
surface and $E$ is a rank two bundle on $S$, of the previous type.

There is a natural rational map
$\Phi:{\cal H} --\ra |{\cal O}_{{\Bbb P}^3}(3)|$,
where ${\cal H}$
is the component of the Hilbert scheme studied in Section 3. The map $\Phi$
sends a scroll $X$ to the base $S$ of the scroll. If $X$ is smooth then $S$ is
its image via the adjunction map. Moreover we have the
following commutative diagram of rational maps:

\

$\put(100,7){${\cal S}_3$}
%\put(114,10){\vector(1,0){40}}
\put(129,10){\line(1,0){10}}
\put(144,10){\vector(1,0){10}}
\put(114,10){\line(1,0){10}}
%\put(114,18){\vector(3,2){40}}
\put(114,18){\line(3,2){10}}
\put(129,28){\line(3,2){10}}
\put(144,38){\vector(3,2){10}}
\put(120,-2){$pf$}
\put(158,7) {$|{\cal O}_{{\Bbb P}^3}(3)|$}
%\put(170,39){\vector(0,-1){20}}
\put(170,39){\line(0,-1){8}}
\put(170,27){\vector(0,-1){8}}
\put(164,45){${\cal H}$}
\put(172,27){$\Phi$}
%\put(155,45){\vector(-3,-2){35}}
\put(125,37){$\Psi$}$

\

\noindent where $\Psi$ is the map which sends the matrix $M=Ax+By+Cz+Dt$ to
the degeneracy locus $X$ as in Section 3.

Let $X$ be a Palatini scroll and let $\Delta$ be the web associated to it.
The fibre of $\Psi$ over $X$
 is $16$-dimensional: its elements correspond to the different choices of a
 base of $\Delta$.
Note that $\Psi$ does not induce any map from the quotient
${\cal S}_{3}/GL(6)$ to ${\cal H}$.

 As a consequence of the interpretation of ${\cal S}_{3}/GL(6)$ seen above,
 $\overline{pf}$ can be
interpreted as a forgetful map and we get a factorization of
$\Phi$ through ${\cal S}_3/GL(6)$, as $\Phi= \overline{pf}\circ
\overline{\Phi}$.

The new map $\overline{\Phi}:{\cal H} --\ra {\cal S}_3/GL(6)$
sends a Palatini scroll $X={\Bbb P}(E)$ to the corresponding pair
$(S,E)$. Since $\dim {\cal H}=44$ and $\dim {\cal S}_3/GL(6)=24$,
the fibres of $\overline{\Phi}$ are $20$-dimensional.

To explain this number, let us observe that the points of the
fibre $\overline{\Phi}^{-1}(S,E)$ are nothing more than the
different embeddings of the projective bundle ${\Bbb P}(E)$ on $S$
in ${\Bbb P}^5$, i.e. the automorphisms of ${\Bbb P}^5$ preserving
$X$.

 We claim that if an automorphism of ${\Bbb P}^5$ preserves
$X$, then it belongs to the subgroup $\overline P$ of $PGL(5)$ of
automorphisms inducing the identity on ${\Bbb P}^3={\Bbb
P}(H^0({\cal O}_S(1)))$.

Indeed, if $G$ is an automorphism of ${\Bbb P}^5$ which preserves
$X$, it preserves also ${\Bbb P}^3$, because $${\Bbb P}^3={\Bbb
P}(H^0({\cal O}_X(K+2H)))\simeq {\Bbb P}(H^4({\cal I}_X(-2))),$$
so it induces on ${\Bbb P}^3$ an action which preserves $S$. But
every automorphism of $S$ is trivial, being $S$ general, so the
restriction of $G$ to ${\Bbb P}^3$ is the identity (see
\cite{ko}). This shows that the fibre of $\overline \Phi$ is
isomorphic to the subgroup $\overline P$, whose dimension is $20$.

%%%%%%%%%%%%%%%%%%% References %%%%%%%%%%%%%%%%%%%%%%%%%%%%%%%%%%%%%%%%%%%
%%%%%%%%%%%%%%%%%%%%%%%%%%%%%%%%%%%%%%%%%%%%%%%%%%%%%%%%%%%%%%%%%%%%%%%%%%

\enddocument